\documentclass[12pt,leqno]{article}

\usepackage[leqno]{amsmath}
\usepackage[normalem]{ulem}
\usepackage{paralist,amsthm,amsfonts,amssymb,graphics, topcapt, amsmath}%
\usepackage{t1enc, pst-poly, pstricks,color,appendix}
\usepackage[active]{srcltx} % SRC Specials: DVI [Inverse] Search
\usepackage{graphicx}

\usepackage{hyperref}% para pdf con hyperlinks
\usepackage{amsmath}
\usepackage{amsfonts}
\usepackage{amssymb}

\divide
\oddsidemargin by 2 
\makeatletter
\newcommand{\address}[2][]{%
  \ifx\@add@ress\@undefined\gdef\@add@ress{\par\par\bigskip}\AtEndDocument{\@add@ress}\fi
  \g@addto@macro\@add@ress{\bigskip\noindent{\small\scshape%
      \ifx#1\empty\else{\bfseries Address of #1:}\ \fi#2}\par\par}}

\renewenvironment{abstract}{\small\quotation\noindent
  {\bfseries \abstractname}}{\endquotation \par}
\newcommand{\footnotetextplain}[1]{\begingroup\def\@thefnmark{}%
  \long\def\@makefntext##1{\parindent 0pt\noindent ##1}\@footnotetext{#1}
  \endgroup}

\newcommand{\AMSsubjclass}[2]{\footnotetextplain{2010
   \emph{Mathematics Subject Classification:} Primary #1, Secondary #2.}}

\newcommand{\keywords}[1]{\footnotetextplain{\emph{Key words and phrases:} #1.}}

\makeatletter \xdef\qedbuit{\qed}

\newcommand{\TeoremaAmbFinalMarcat}[1]{%
  \expandafter\gdef\csname end#1\endcsname{\qedbuit\@endtheorem}}

\newtheorem{theo}{Theorem}[section]

\theoremstyle{definition}
 \TeoremaAmbFinalMarcat{defi}
 \TeoremaAmbFinalMarcat{ex}
\newtheorem{rem}[theo]{Remark} \TeoremaAmbFinalMarcat{rem}
 \TeoremaAmbFinalMarcat{hrem}
 \TeoremaAmbFinalMarcat{conv}
 \TeoremaAmbFinalMarcat{prob}
\newtheorem{conj}[theo]{Conjecture} \TeoremaAmbFinalMarcat{conj}
 \TeoremaAmbFinalMarcat{conj}

\newenvironment{proclama}[1]{
                \par\vspace{\topsep}\noindent{\bf #1}
                \begin{em}}
                {\end{em}\par\vspace{\topsep}}

\newcommand{\start}[2]{\begin{#1}\label{#2}}

\newcommand{\theoc}[1]{Theorem~\ref{#1}}
\newcommand{\propc}[1]{Proposition~\ref{#1}}

\newcommand{\lemc}[1]{Lemma~\ref{#1}}

\newcommand{\figc}[1]{Figure~\ref{#1}}

\newcommand{\defc}[1]{Definition~\ref{#1}}

\def\@enum@{\list{\csname label\@enumctr\endcsname}%
           {\usecounter{\@enumctr}\def\makelabel##1{\hss\llap{##1}}
           \itemsep=2pt\parsep=0pt\topsep=3pt plus 1pt minus 1 pt}}

\newenvironment{numlist}{\enumerate[(1)]}{\endenumerate}

\def\map#1#2#3{\mbox{${#1}\colon {#2} \longrightarrow {#3}$}}
\def\Smap#1#2{\mbox{${#1}\colon{#2} \longrightarrow {#2}$}}
\def\u{\mathsf{U}}
\def\v{\mathsf{V}}
\def\w{\mathsf{W}}
\def\x{\mathsf{X}}
\def\y{\mathsf{Y}}

\def\p{\mathsf{P}}
\def\q{\mathsf{Q}}

\def\UU{\mathcal{U}}
\def\VV{\mathcal{V}}
\def\WW{\mathcal{W}}

\def\V{\mathbb{V}}
\newcommand{\vvv}{\pi^\ast}
\def\W{\mathbb{W}}
\def\bb{\mathbf{LP}}

\def\ov#1{\overline{#1}}

\def\wh#1{\widehat{#1}}
\def\wv{\widehat{\v}}

\def\id{\mathop\mathrm{Id}}
\def\sign{\mathrm{sign}}

\def\cob{\delta}
\def\cib{\Delta}

\def\om{\omega}
\def\sss{\mathsf{s}}

\newcommand{\aq}{\ensuremath{\mathbb{A}}_q}

\newcommand{\Z}{\ensuremath{\mathbb{Z}}}

\newcommand{\m}{\mathsf{M}}
\def\lu{\ell_\mathsf{U}}
\def\lv{\ell_\mathsf{V}}
\def\lw{\ell_\mathsf{W}}

\newcommand{\xx}{
 \rput(-0.5,1.3){$\w$}
 \rput(-0.5,0.3){$\u$}
 \rput(0,1.5){
 \rput(1.25,0){$\v_k$}
 \rput(3,0){$\v_k$}
 \rput(5,0){$\v_k$}
\rput(6.75,0){$\v_{\scriptscriptstyle k,h}$}}
\psline(0,2)(7.5,2)
\psline[linecolor=gray,linewidth=3pt](0,2)(0,1)
\psline(2,1)(2,2)
\psline(4,1)(4,2)
\psline[linewidth=3pt](6,1)(6,2)
\psline[linecolor=gray,linewidth=3pt](7.5,1)(7.5,2)
\psline(0,1)(7.5,1)

 \rput(1.3,0){
  \psline[linecolor=gray,linewidth=3pt](0.1,1)(0.1,0)
 \rput(1,0.3){$\v_i$}
 \rput(3,0.3){$\v_i$}
 \rput(5,0.3){$\v_i$}
  \rput(7.2,0.3){$\v_{\scriptscriptstyle i\!,\!j}$}
  \psline[linecolor=gray,linewidth=3pt](7.5,0)(7.5,1)%
\psline(7.5,1)(0.2,1)
\psline(2,0)(2,1)
\psline(4,0)(4,1)
\psline[linewidth=3pt](6,0)(6,1)
\psline(0.1,0)(7.5,0)
}

\psdot(0.2,0.7)
\psline{->}(0.2,1.2)(0.2,0.8)

\psdot[dotstyle=o](0.2,1.2)
\rput(0.25,1.45){$v_{k}$}

\psdot(5.6,0.7)
\psline{->}(5.6,0.8)(5.6,1.2)
\psdot(5.6,1.2)

\rput(5.6,1.4){$v_h$}
\psline{->}(7.5,0.8)(5.7,0.8)

\psdot[dotstyle=o](7.6,1.2)
\rput(7.8,0.4){$u_{\scriptscriptstyle \lu\!-\!r}$}

\psline{->}(7.6,1.1)(7.6,0.8)
\rput(6.6,0.6){${\scriptstyle \lv}$}

\psdot(7.6,0.7)
\rput(7.7,1.45){$v_k$}

  \psline[linecolor=gray]{[-}(0.1,0.4)(0.4,0.4) \rput(0.6,0.4){$r$}\psline[linecolor=gray]{-]}(0.7,0.4)(1.2,0.4)
}

\title{Algebraic characterization of simple closed curves via Turaev's cobracket}

\author{Moira Chas and Fabiana Krongold}
\date{ }
\makeindex

\begin{document}
\maketitle
%\tableofcontents
\begin{abstract}  The vector space $\V$ generated by the conjugacy classes in the fundamental group  of an orientable surface has a natural Lie cobracket $\map{\delta}{\V}{\V\otimes \V}$. For negatively curved surfaces, $\delta$ can be computed from a geodesic representative as a sum over  transversal self-intersection points. In particular $\delta$ is zero for any power of an embedded
simple closed curve. Denote by Turaev(k)  the statement that $\delta(x^k) = 0$ if and only if the non-power conjugacy class $x$ is represented by an embedded curve.
Computer implementation of the cobracket $\delta$ unearthed counterexamples to Turaev(1) on every surface with negative Euler characteristic except genus zero surfaces.
Computer search has verified Turaev(2) for hundreds of millions of the shortest classes. In this paper we prove Turaev(k) for $k=3,4,5,\dots$ for surfaces with boundary.
Turaev himself introduced the cobracket in the 80's and wondered about the relation with embedded curves, in particular asking if a statement equivalent to Turaev (1) might be true.

We give an application of our result to the curve complex. We show that a permutation of the set of free homotopy classes that commutes with the cobracket and the power operations is induced by an element of the mapping class group.

\end{abstract}
\AMSsubjclass{57M05}{17B62}

\keywords{surfaces, conjugacy classes, Lie bialgebras, self-intersection number, embedded curves}

\thanks{Supported by NSF grant DMS  1105772}

\section{Introduction}

Turaev \cite{T}, defined a Lie coalgebra structure in the vector space generated by free homotopy classes of non-trivial closed curves on an orientable surface. The cobracket of a closed curve is a sum  of certain terms over self-intersection points of the curve. Therefore, the cobracket of a simple closed curve is zero. It is not hard to see that the cobracket of a power of a simple closed curve is zero. Thus, he asked whether the converse is true. Le Donne \cite{ledonne} answered this question positively for the case of surfaces of genus zero. The cobracket is zero in the torus, where all the curves are power of a simple one.  For all other surfaces of positive genus, the answer to this question is negative (see \cite{chas}). The main goal of our work here is to show that one may reformulate Turaev's question so that it has a positive answer, at least for surfaces with boundary. Moreover, it also can be shown in that case that the cobracket of a power larger than two of a curve "counts" the self-intersection number of its class. (Recall that the \emph{self-intersection number of a free homotopy class $\alpha$} is the smallest number of self-intersection points of a representative of $\alpha$ whose intersection points are all transversal double points.) In particular, Turaev's cobracket yields a characterization of simple closed curves. More precisely, we prove:

\begin{proclama}{Main Theorem} Let $\VV$ be a free homotopy class of curves on an oriented surface with boundary.
 Then  $\VV$ contains a power of a simple curve if and only if the Turaev cobracket of $\VV^3$ is zero. Moreover, if $\VV$ is a non-power and $p$ is an integer larger than three then the number of terms of the cobracket of $\VV^p$ (counted with multiplicity) equals $2p^2$ times the self-intersection number of $\VV$.\end{proclama}

Our main tool is the combinatorial presentation of Turaev's cobracket in \cite{chas}. Using this presentation, we list the terms of the cobracket of powers of a cyclic word and we can show that if the power is large enough, these terms do not cancel.

We also give an application of our result to the curve complex.

The problem of characterizing free homotopy classes containing simple curves has a rich history, starting probably in 1895 with Poincar\'e \cite{poincare}, who showed   that a free homotopy class contains a simple representative if and only if the unique smooth geodesic representative of the class is  simple, and that a geodesic is simple if and only if  each pair of liftings of the geodesics are disjoint. A  purely group-theoretical method of determining simplicity  was given in the 1960's by Zieschang \cite{z1,z2}. Chillingworth \cite{chi1,chi2} in 1972 gave a geometric algorithm in terms of the winding number.  In 1983, Turaev and Viro \cite{tv} gave an  algebraic method to count minimal intersection and self-intersection using the fundamental group.  In 1984, Birman and Series \cite{BS} gave an algorithmic characterization of simple closed curves in surfaces with non-empty boundary. Also in 1984 Cohen and Lusig \cite{cl},  extended Birman and Series ideas to count the minimal number of intersection and self-intersections of curves on such surfaces and Lustig \cite{l} gave a refinement of this  algorithm for the case of closed surfaces. Hass and Scott \cite{hs1} gave a geometric algorithm to produce a minimal representative of a free homotopy class. Arettines \cite{Ar} gave a combinatorial algorithm for the same purpose.  We  used the Goldman bracket to give a characterization of simple class on surfaces with boundary \cite{CK}. Recently, in 2010 Cahn \cite{Cahn} gave another characterization by generalizing Turaev's cobracket  to garlands.

Much of the interest in characterizing simple conjugacy classes was stimulated by the Jaco-Stallings group theory statement equivalent to the Poincar\'e conjecture. Our bracket and cobracket characterizations led to the generalization of  Goldman-Turaev   on surfaces to the String Topology Lie bialgebra for higher dimensional manifolds. Of course, it is still an interesting problem to prove the Jaco-Stallings \cite{jaco, stallings} criterion directly in terms of geometric group theory.

This paper is structured as follows: We start by describing in Section \ref{tur} the vector space of cyclic reduced words  and  the combinatorial cobracket. In Section \ref{not conjugate}, using techniques similar to those   in \cite{CK}, we prove that certain pairs of cyclic
words cannot be conjugate (\propc{main}).
This result distinguishes the conjugacy classes of various sets of linear words. In Section~\ref{turaev}, we use Proposition~\ref{main} to prove our Main Theorem. In Section~\ref{ccc} we give an application of our Main Theorem to the curve complex.
In Appendix~\ref{curves},
 we recall the geometric definition of Turaev's cobracket.
 In Appendix \ref{review c}, we give a geometric interpretation of the definitions of Section~\ref{tur} (Appendix~\ref{review c} is not necessary for the logic  but it suggests the geometric picture of what is being rigorously proved by the combinatorics). Finally, in Appendix \ref{cc} we review the definition of the curve complex.

\textbf{Acknowledgment}
We are very thankful to the referee for the careful reading of our draft and the many useful suggestions which have  made this paper much more readable.

\section{The Turaev Lie coalgebra on the vector space generated by cyclic reduced words}\label{tur}
The results of this section are purely combinatorial. See Appendix~\ref{review c} for a geometric interpretation.

\subsection{Definitions and preliminary results}\label{Notation}

In this section, we will state some results from \cite{chas} which are used throughout this work. For convenience, we reformulate these ideas in terms of the notation of \cite{CK}.

For each positive integer $q$,  a \emph{$q$-alphabet $\aq$}, or, briefly,
an \emph{alphabet}, is the set of  $2q$ symbols, called \emph{letters} $\{a_1,a_2,
\dots,a_q,\ov{a}_1, \ov{a}_2, \dots, \overline{a}_q\}$, endowed with
a fixed linear order.  For each letter $v$,
$\ov{\ov{v}}=v$.

A \emph{linear word in $\aq$}  is a  finite sequence of
symbols $v_0v_1\dots v_{n-1}$ such that $v_i$ belongs to $\aq$ for
each $i \in \{0,1,\dots,n-1\}.$  The empty word is a linear word with zero letters. Let $\v=v_0v_1\dots v_{n-1}$ be a linear word. By
definition, $\ov{\v}=\ov{v}_{n-1} \ov{v}_{n-2}\dots \ov{v}_0$.   The
linear word $\v$ is \emph{freely reduced} if $v_i \ne \ov {v}_{i+1}$
for each $i \in \{0,1,\dots,n-1\}$. If $\v$ is freely reduced and
$v_{n-1}\ne \ov{v}_0$ then $\v$ is \emph{cyclically reduced}.
Consider the equivalence relation on the set of linear words,
generated by the pairs of the form $(\v,\w)$ such that $\v$ is a
cyclic permutation of $\w$ or $\v=\w v \overline{v}$ where $v$ is a
letter in $\aq$. The equivalence classes  under this equivalence
relation  are called \emph{cyclic words}. (Observe that these are
the conjugacy classes of the free group generated by $a_1, a_2,
\dots, a_{q}$). If $\v$
is a (not necessarily reduced) linear word, we denote the
equivalence class of $\v$ by $\widehat{\v}$. A \emph{ linear representative} of $\widehat{\v}$ is a cyclically
reduced linear word in $\wv$.
A cyclic word $\widehat{\v}$ is  \emph{nonpower} if it is not  a proper power of another word.
A linear subword of ${\widehat\v}$ is a linear subword of a linear representative of ${\widehat\v}$.
The \emph{length of $\v$}, denoted by $\ell_\v$ or $\ell(\v)$  is  the number of symbols $\v$ contains. The \emph{length of  the cyclic word $\widehat{\v }$} is  the length of  any linear representative of $\widehat{\v }$.  (Note that all of them have the same length) Thus, if $\v$ is cyclically
reduced then the length of $\widehat{\v }$ equals $\ell_\v$.

 When dealing with letters denoting linear
words $\v=v_0v_1\dots v_{n-1}$ , subindices of
letters will be considered mod the length of $\v$, which is $n$.

\start{defi}{oriented}
Let $r$ be an integer greater than three and $(x_0$, $x_1$,\ldots, $x_{r-1})$ a sequence of $r$ letters in the ordered alphabet $\aq$. The sequence is \emph{ positively} (resp. \emph{negatively})
\emph{ oriented} if for some a cyclic permutation  the sequence is strictly increasing (resp. strictly decreasing). Note that this implies that there are no consecutive repetitions in the sequence.
The sign of  $(x_0,x_1,\ldots, x_{r-1})$ , denoted by $\sign(x_0,x_1,\ldots, x_{r-1})$, is 1 (resp. -1) if  $(x_0,x_1,\ldots, x_{r-1})$ is positively (resp. negatively) oriented and zero otherwise.
\end{defi}

From now on,   fix $\v=v_0v_1\dots v_{n-1}$, a cyclically reduced linear word of positive length $n$. Let $\p$ and $\q$ be linear subwords of $\widehat{\v}$ of equal length possibly overlapping.

\subsection{Linked pairs}

In this section  the definition of linked pairs is recalled. Examples of such pairs can be found in Tables~\ref{table} and \ref{table1},  Appendix~\ref{review c}.

\start{defi}{linked} \cite[Definition 2.1]{chas}  The ordered pair $(\p,\q)$  is a {\em linked pair} if there exists non-negative integers $i$ and $j$ smaller that $n$ such that one of the following holds:
\begin{numlist}
\item $\p=v_{i-1}v_i$, $\q=v_{j-1}v_j$  and  $(\ov{v}_{i-1},\ov{v}_{j-1},v_i,v_j)$ is oriented. (Note that in particular $v_{i-1} \ne v_{j-1}$, $v_{i} \ne v_{j}$)
\item $\p=v_{i-r-1}\y  v_{i}$, $\q=v_{j-r-1}\y  v_{j}$ for some $r$ in  $\{1,2,\dots,n-2\}$. Moreover,
$\y=v_{i-r}v_{i-r+1}\cdots v_{i-1}=v_{j-r}v_{j-r+1}\cdots v_{j-1}$   and the sequences $(\ov{v}_{i-r-1} ,\ov{v}_{j-r-1} ,v_{i-r})$ and $(v_{i}, v_{j}
,\ov{v}_{i-1})$ have the same orientation, both positive or both negative.

\item $\p=v_{i-r-1}\y  v_{i}$, $\q=v_{j-1}\overline{\y}  v_{j+r}$ for some $r$ in  $\{1,2,\dots,n-2\}$. Moreover, $\y=v_{i-r}v_{i-r+1}\cdots v_{i-1}=\ov{v}_{j+r-1}\ov{v}_{j+r-2}\cdots \ov{v}_{j}$  and the sequences $(v_{j+r}, \ov{v}_{i-r-1},v_{i-r})$ and
$(\overline{v}_{i-1},\ov{v}_{j-1},v_i)$ have the same
orientation, , both positive or both negative.
\end{numlist}

The \emph{sign} of a linked pair  ${(\p,\q)}$ is given by the formula
$\sign(\p,\q)= \sign(\ov{v}_{i-1} ,v_i ,v_j) $.

The set of linked pairs of a cyclic word $\widehat{\v}$ is denoted by $\bb(\widehat{\v})$.

\end{defi}

The next statement follows straightforwardly from  \defc{linked}.

\start{lem}{pair}
For each $({  \p, \q})$in  $\bb(\wv)$ there exist two integers  $i$ and $j$  in $\{0,1,\ldots,n-1\}$ such that the following holds:

\begin{numlist}
\item If $({  \p, \q})$ is a linked pair as in
\defc{linked} (1) or (2) then $\p=v_{i-r-1}v_{i-r}\dots v_{i}$ and
$\q=v_{j-r-1}v_{j-r}\dots v_{j}$, for some  integer $r \in
\{0,1,\ldots,n-2\}$.
\item If $({  \p, \q})$ is a linked pair as in \defc{linked} (3) then $\p=v_{i-r-1}v_{i-r}\dots v_{i}$ and
$\q=v_{j-1}v_j\dots v_{j+r}$,  for some  integer $r \in
\{1,2,\ldots, n-2\}$.
\end{numlist}
Moreover,  the map  $\phi: \bb(\wv)\to\{0,1,\ldots, n-1\} \times \{0,1,\ldots, n-1\}$ defined by $({
\p, \q})\mapsto (i,j)$ is injective.
\end{lem}

Since $\v$ is fixed, we can (and will)  identify each linked pair $(\p,\q)$ with the pair of integers $\phi(\p,\q)$ in \lemc{pair}.
Say that $(i,j)$ is a \emph{linked pair of type (1) or (2)} (respectively, \emph{(3)}) if $(i,j)=\phi(\p,\q)$ for some $(\p,\q)$ in $\bb(\wv)$ satisfying \defc{linked}(1) or (2) (respectively (3)).

\start{rem}{choice} The map $\phi$ of \lemc{pair} is defined to reflect a particular pair of segments where the representative of the class  intersects. The choice of the pair of segments is somewhat arbitrary for linked pairs of type (2) and (3). In case (2), the words $(\p,\q)$ correspond to two arcs of the representative  made up of $r+1$ segments each. Our definition of $\phi$ corresponds to having the intersection point in the last pair of segments of the arcs, that is, in the intersection of the arc from $v_{i-1}$ to $v_i$ and from $v_{j-1}$ to $v_j$.
\end{rem}

\start{lem}{extremo}If $(i,j)$ is a linked pair of  $\wv$ then $v_{i}\ne v_{j}$ and
$\ov{v}_i\ne v_{j-1}$. In particular $i \ne j$.
\end{lem}

\begin{proof}
If $(i,j)$ is a linked pair of type (1) then  $(\overline{v}_{i-1},\overline{v}_{j-1},v_i,v_j)$  is positively
or negatively oriented. Hence,  all consecutive letters must be different: $v_{i-1}\ne v_{j-1}$, $v_i\ne v_j$ and $\overline{v}_{j-1}\ne v_i$ .

If $(i,j)$ is a linked pair of type (2), then  $(v_i,v_j,v_{i-1})$ is positively or negatively oriented. This implies that $v_i \ne v_j$ and since $v_i \ne \overline{v}_{i-1}$ because $\wv$ is cyclically reduced and $v_{i-1}=v_{j-1}$ by hypothesis, the result follows.

Finally if $(i,j)$ is of type (3), the last letter of $\y$ (see Definition \ref{linked}) is $v_{i-1}$ and the first letter of $\overline{\y}$ is $v_j$. Thus, $v_{i-1}=\overline{v}_j$. On the other hand, $(\overline{v}_{i-1}, v_i,\overline{v}_{j-1})$ is positively or negatively oriented.

\end{proof}

For each pair $i$ and $j$  in $\{0,1,\dots, n-1\}$ denote
\[\v_{i,j}=\left\{\begin{array}{ll}v_{i} v_{i+1}\ldots v_{j-1} & \hbox{ if $i<j$, and}\\
                             v_{i}  v_{i+1}\ldots v_{n-1}v_{0}v_1\ldots v_{j-1} & \hbox{ if $i\ge j$.}\end{array}
            \right.\]
If $i=j$, $\v_i$ is used instead of $\v_{i,i}$.  Note that in this case $\v_i$ is just the linear word associated to $\v$ starting at position $i$.

The next result is a direct consequence of \defc{linked}.
\start{lem}{symmetry} \begin{numlist} \item If $(\p,\q)$ is a linked pair then $(\q,\p)$ is a linked pair.
However, $(i,j)$ and $(j,i)$ are both linked pairs  if an only if $(i,j)$ is a linked pair of type (1) or (2).
\item Let $(i,j)$ be a linked pair.  The linear word if $\v_{j,i}$ is cyclically reduced if an only if $(i,j)$ is a linked pair of type (1) or (2).
\item If $(i,j)$ is  a linked pair of type (3), then there exists  $r$ in $\{1,2,\dots,\ell_\v-1\}$ such that $r$ is the largest integer verifying  $v_{i-r}\dots v_{i-1}=\overline{v_{j}\dots v_{j+r-1}}$ and $(j+r,i-r)$ is a linked pair.
\end{numlist}
\end{lem}

\subsection{The ``combinatorial'' Turaev cobracket}

 In this section $\Delta$ stands for Turaev's topological cobracket (see Appendix \ref{curves}).

\start{defi}{cob12}
To each linked pair $(\p,\q)$ in $\bb(\widehat{\v})$ we associate two cyclic words,
$\cob_1(\p,\q)=\wh{\v}_{\hbox{\tiny $\!\!i\!,\!j$}}$ and $\cob_2(\p,\q)= \wh{\v}_{\hbox{\tiny $\!\!j\!,\!i$}}$.
\end{defi}

Denote by  $\V$  the vector space generated by the set of non-empty cyclic words in
$\aq$.
The class of the empty word will be identified with the zero element in $\V$.

\start{defi}{cobr}
Let $\map{\delta}{\V}{\V\otimes \V}$  be the  map defined by
\begin{equation*}
    \cob(\wv)=\sum_{(\p,\q)\in \bb(\wv)}\sign(\p,\q)\,\,\cob_1(\p,\q)\otimes
    \cob_2(\p,\q),
\end{equation*}
for every $\wv$ in $\V$ and extended by linearity.
\end{defi}

Recall that if $\Sigma$ is a surface with non-empty boundary, the set of non-trivial, free homotopy classes of closed, oriented curves on $\Sigma$ is in bijective correspondence with the set of non-empty,  cyclic reduced words on a set of free generators of the fundamental group of $\Sigma$ and their inverses.

\start{theo}{iso} (\cite[Proposition 4.1]{chas}) If $\Sigma$ is an
orientable surface with non-empty boundary then $\delta(\wv)=\Delta(\wv)$ for each $\wv$ in $\V$.
\end{theo}

The next theorem was stated as  Remark 3.10 in \cite{chas}.

\start{theo}{intersections} If $\VV$ is a
nonpower cyclic word and $\v$ a linear representative of $\VV$ then  the number of linked pairs of $\wv$ equals twice the self-intersection number of $\VV$.
\end{theo}
\begin{proof}
By \cite[Proposition 3.3]{chas}, there exists a closed curve $\alpha$ representing $\VV$ with no singular $1$-gon and no singular bigons. By \cite[Theorem 4.2]{hs} the number of self-intersection points of $\alpha$ equals the self-intersection number of $\VV$.
On the other hand, by \cite[Theorem 3.9]{chas} the self-intersection points of $\alpha$ are in bijective correspondence with the set of pairs of linked pairs of $\VV$ of the form $\{(\p,\q),(\q,\p)\}$. Finally, if $(\p,\q)$ is a linked pair then $\p \ne \q$. Thus the desired conclusion follows.
 \end{proof}

The next result which establishes the exact relation between the linked pairs of a nonpower cyclic word $\widehat{\v}$ and its power $\widehat{\v}^p$, follows directly from \defc{linked}.

\start{lem}{linkedvp}
Let $\v$ be a cyclically reduced word of length $n$ and $p$ be a positive integer. Then,
$$\bb(\wv^p)=\Big\{\big(i+tn,j+sn\big): (i,j)\in \bb(\wv),  t,s\in\{0,1,2,\dots,p-1\}\Big\}.$$
%In particular   \(|\bb(\wv^p)|=p^2 |\bb(\wv)|\).
\end{lem}

Denote by  $\V$  the vector space generated by the set of non-empty cyclic words in $\aq$.
The class of the empty word will be identified with the zero element in $\V$.

\start{prop}{delta} If $\v$ is cyclically reduced linear word and $p$ is a positive integer, then the map
$\map{\delta}{\V}{\V\otimes \V}$ of \defc{cobr} satisfies the following equality.
\begin{equation}{\label{eq delta}}
\delta(\widehat{\v}^p)= p \cdot\hspace{-5mm}\sum_{ \substack{(i,j)\in
                \bb(\wv) \\ 0\le s\le p-1}}\hspace{-5mm} \sign(i,j) \hspace{0.1cm} \wh{\Big(\v_i^s\v_{i,j}\Big)}\otimes
                 \wh{\Big(\v_j^{p-s-1}\v_{j,i} \Big)}.
\end{equation}
\end{prop}
\begin{proof}  Take $i, j \in \{0,1,\dots,n-1\}$ and $t, s \in \{0,1, \dots, p-1\}$. By \lemc{extremo}, $i \ne j$. It is not hard to see that

\[ \v^p_{i+tn,j+sn}=\left\{\begin{array}{ll}
				\v_i^{s-t}\ \v_{i,j} & \hbox{ if $i<j$ and $t \le s$,}\\
		         	\v_i^{p-t+s}\ \v_{i,j} & \hbox{ if $i<j$ and $t>s$,}\\	
                                 \v_i^{s-t-1}\ \v_{i,j} & \hbox{ if $i>j$ and $t<s$,}\\
                                  \v_i^{p-t+s-1}\ \v_{i,j} & \hbox{ if $i>j$ and $t \ge s$,}\\
                             	\end{array}
            \right.\]
The desired result follows from  \lemc{linkedvp}.
\end{proof}

\section{Certain words are not conjugate}\label{not conjugate}

In this subsection we prove combinatorially that certain pairs of cyclic words we construct out of  $\v$  cannot be
conjugate. Since the conjugacy classes of these words appear in the
terms of the cobracket $\delta(\widehat{\v}^p)$,
the fact that these words are not conjugate implies that these terms cannot cancel one
another.

 \start{prop}{main} Let $\v$  be a  cyclically reduced linear word  such that $\widehat{\v}$ is   nonpower and  $(i,j)$ and $(k,h)$ two distinct linked pairs of $\v$. If
$m$ is a positive integer greater than or equal to two, then, for
any positive integer $l$,
$$
 \widehat{\v_i^{m}\v}_{i,j}\not= \widehat{\v_k^{l}\v}_{k,h}
%\ \  \hbox{ and }\ \ \widehat{\v_{j,i}\v_i^{m}}\not= \widehat{\v_{h,k}\v_k^{r}}.
$$
\end{prop}

\begin{proof}

Let  $\u$  and $\w$ denote the linear words $\v_i^m\v_{i,j}$ and
$\v_k^l\v_{k,h}$ respectively.  By hypothesis, $\v$ is a cyclically
reduced linear word, therefore $\u$ and $\w$ are reduced. Moreover,
by \lemc{extremo}, $v_i\not=\overline{v}_{j-1}$ and
$v_k\not=\overline{v}_{h-1}$. Hence $\u$ and $\w$ are cyclically
reduced and
\begin{align} \label{d1}
m\cdot\lv< \lu <(m+1)\cdot \lv,\\
\nonumber l\cdot \lv< \lw <(l+1)\cdot \lv.
\end{align}

Set $\u=u_0u_1\ldots u_{\lu-1}$ and $\w=w_0w_1\ldots w_{\lw-1}$. The next equalities follow straightforwardly.
\begin{align}\label{p1}
 u_t&=u_{t+\lv} ,  \hbox{   for all $t$ in  $\{0, 1, \ldots , \lu -\lv-1\}$,}\\
\label{p2} w_t&=w_{t+\lv} ,  \hbox{   for all $t$ in $\{0, 1, \ldots , \lw -\lv-1\}$,}\\
\label{u} u_{\lu-\lv}&= v_{i+\lu-\lv}=v_j, \\
\label{w} w_{\lw-\lv}&=v_{k+\lw-\lv}=v_h.
\end{align}

We argue by contradiction: Assume that $\widehat\u= \widehat\w$.
Since the linear words $\u$ and $\w$ are linear representatives of
the same cyclic word, they have equal length, that is, $\lu=\lw$.
Furthermore, there exists $r$ in $\{0,1,\ldots,\lu-1\}$ such that
$\u=\w_{r}$. The proof is split into three cases:

\begin{enumerate}

\item $r=0$. Here, $\u=\w$. Hence the initial and final subwords of length
$\lv$ of $\u$ and $\w$ are equal i.e., $\v_i=\v_k$ and
 $\v_j=\v_h$. Distinct cyclic permutations of a nonpower cyclic reduced word are distinct, so $i=k$ and $j=h$, contradicting the hypothesis $(i,j)\ne (k,h)$.
\item    $0<r\le \lu-\lv$. By Equations (\ref{p2}) and (\ref{w}), since  $m$ is positive,
\[v_{k}=w_0=w_{\lw}=u_{\lw-r}= u_{\lw-r-\lv}=w_{\lw-\lv}=v_h,\]
%\[v_{j-1}=u_{-1}=w_{-1 + r}=w_{\lv+ r-1}= u_{\lv-1}=v_{i-1},\]
which contradicts \lemc{extremo}. See \figc{Case 2}.

\begin{figure}[htbp]
 \begin{pspicture}(14,2)
 \rput(2,0){\xx}
   \end{pspicture}
  \caption{
   Case  $0<r\le \lu-\lv$ (here, $m=l=3$.)}
  \label{Case 2}
 \end{figure}
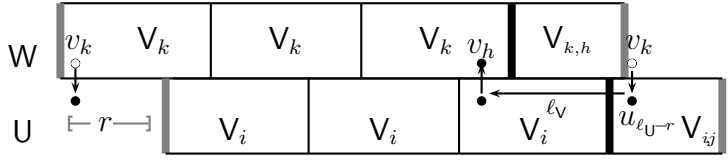

\item    $\lu-\lv< r< \lu$. Since $m \ge 2$, $\lv-1<\lu-\lv$. Hence $0\le \lu-r -1< \lu-\lv$.
By Equations (\ref{p1}) and (\ref{u}),
%\[v_{k-1}=w_{\lw+\lv-1}=u_{\lu+\lv-r-1}=u_{\lu-r-1}=w_{\lw-1}=v_{h-1},
%\]
\[v_{j}=u_{\lu-\lv}=w_{\lw-\lv+r}=w_{\lw+r}=u_{\lu}=u_0=v_i,
\]
contradicting \lemc{extremo}.

\end{enumerate}
\end{proof}
\begin{rem}\label{ejemplo} The lower bound $2$ for $m$ in  \propc{main} is sharp.
 Indeed, consider the ordered alphabet $\{a,b,\overline{b},\overline{a}\}$ and the word $\v=abaabab$. The  pairs $(0,1)$ and
$(5,6)$ are linked. Nevertheless, the subword of $\v^2$ given by
$\v_0\v_{0,1}=abaababa$ is a cyclic permutation of $\v_5\v_{5,6}=ababaaba $. (Observe that the corresponding terms in the cobracket of $\v^2$,  $\widehat{\v_0\v}_{0,1}\otimes
\widehat{\v}_{1,0}$ and $\widehat{\v_5\v}_{5,6}\otimes \widehat{\v}_{6,5}$ are different because
$\v_{1,0}=baabab$ and $\v_{6,5}= babaab$ which implies
$\widehat{\v}_{1,0} \ne \widehat{\v}_{6,5}=babaab$.)
\end{rem}

 The \emph{Manhattan norm of an element $x$
of $\V \otimes \V$} denoted by $\m(x)$
 is the sum of the absolute values of the coefficients of the expression of
 $x$ in the basis of $\V\otimes \V$ consisting in the set of tensor products of pairs of cyclic words
 in $\aq$. Thus if $x=c_1 \VV_1\otimes\WW_1+c_2\VV_2\otimes \WW_2+\cdots+c_l \VV_l\otimes\WW_l$ where for each $i,j \in \{1,2,\dots,l\}$,
 $c_i \in \Z$, $\VV_i$ and $\WW_i$ are  cyclically reduced words and $\VV_i \ne \VV_j$ or $\WW_i \ne \WW_j$ when $i \ne j$ then $\m(x)=|c_1|+|c_2|+
 \cdots+|c_l|$.

\start{prop}{contar} Let $\v$  be a  cyclically reduced linear word and let $p$ an integer larger than three. Then the Manhattan norm of $\delta(\widehat{\v}^p)$ equals $p^2$ times
the number of elements in $\bb(\wv)$.\end{prop}

 \begin{proof} A simple argument shows that if the result holds for nonpower words, then it holds for all words.  We assume then that $\widehat{\v}$ is nonpower.
Let  $(i,j)$ be an element of $\bb(\wv)$. The terms of $\delta(\widehat{\v}^p)$ corresponding to $(i,j)$ in Equation (\ref{eq delta}) are
$\widehat{\v_i^m\v_{i,j}} \otimes \widehat{\v_j^{\hbox{\tiny $p\!\!-\!\!m\!\!-\!\!1$}}\v_{ji}}$ where $m \in \{0,1,2, \dots , p-1\}$.
Since the lengths of the first factors of these elements differ by a non-zero multiple of $n$, these terms are all distinct.

We claim that if $(k,h) \in \bb(\wv)$,  $(k,h)\not=(i,j)$, and $l \in \{0,1,2,\dots, p-1\}$ then
$$\widehat{\v_i^m\v_{i,j}} \otimes \widehat{\v_j^{\hbox{\tiny $p\!\!-\!\!m\!\!-\!\!1$}}\v_{ji}} \ne \widehat{\v_k^l\v_{k,h}} \otimes \widehat{\v_h^{\hbox{\tiny $p\!\!-\!\!l\!\!-\!\!1$}}\v_{hk}}$$
Clearly, this claim implies the desired result.

To prove the claim observe that, since $p\ge 4$, either $m \ge 2$ or $p-m -1 \ge 2$. If $m\ge 2$, the result follows from by \propc{main}.  Hence, we can assume
 that $p-m -1 \ge 2$. If  $\widehat{\v_i^m\v_{i,j}}$ and $\widehat{\v_k^l\v_{k,h}}  $ are distinct the claim follows. Thus we can also assume that $\widehat{\v_i^m\v_{i,j}}=\widehat{\v_k^l\v_{k,h}}$  and therefore $m=l$. We complete the proof by showing that the assumption $\widehat{\v_j^{\hbox{\tiny $p\!\!-\!\!m\!\!-\!\!1$}}\v_{ji}} = \widehat{\v_h^{\hbox{\tiny $p\!\!-\!\!m\!\!-\!\!1$}}\v_{hk}}$
 leads to a contradiction.

Assume first that $(i,j)$ is  a linked pair of type (1) or (2). By \lemc{symmetry}(1)  $(j,i)$ is also a linked pair
of type (1) or (2) and  $(h,k)\not=(j,i)$ (since by hypothesis $(i,j)\not= (k,h)$).
 By  \propc{main}, $\widehat{\v_j^{\hbox{\tiny $p\!\!-\!\!m\!\!-\!\!1$}}\v_{j,i}}  \ne \widehat{\v_h^{\hbox{\tiny $p\!\!-\!\!m\!\!-\!\!1$}}\v_{h,k}} $. Therefore, we can assume that  $(i,j)$ is a linked pair of type (3).
In this case, by  \lemc{symmetry}(2), the word $\v_j^{p-m-1}\v_{ji}$ is not cyclically reduced. Hence,  the sum of the lengths of $\widehat{\v_i^m\v_{i,j}} $ and  $\widehat{\v_j^{p-m-1}\v_{ji}}$ is strictly smaller than  $ p\cdot \ell_v$.

 Therefore, the sum of the lengths of $\widehat{\v_k^l\v_{k,h}}$ and $\widehat{\v_h^{\hbox{\tiny $p\!\!-\!\!l\!\!-\!\!1$}}\v_{hk}}$ is also $p\cdot\ell_\v-2r$.
By \lemc{symmetry}(3),  $(k,h)$ is  a linked pair of type (3) and $(h+r,k-r)$ is linked. Hence,
$$\widehat{\v_j^{\hbox{\tiny $p\!\!-\!\!m\!\!-\!\!1$}}\v_{\hbox{\tiny
$\!\!j\!,\!i$}}}=\widehat{\v_{j+r}^{\hbox{\tiny $p\!\!-\!\!m\!\!-\!\!1$}}\v}_{\hbox{\tiny $j\!\!+\!\!r\!,\!i\!\!-\!\!r$}} =\widehat{\v_h^{\hbox{\tiny $p\!\!-\!\!l\!\!-\!\!1$}}\v_{h,k}}= \widehat{\v_{h+r}^{\hbox{\tiny $p\!\!-\!\!l\!\!-\!\!1$}}\v}_{\hbox{\tiny $h\!\!+\!\!r\!,\!k\!\!-\!\!r$}}.$$
By \propc{main}, $(j+r,i-r)=(h+r,k-r)$ and hence $(i,j)=(k,h)$ contradicting the hypothesis of the claim.

\end{proof}

\start{prop}{simple} For each cyclically reduced linear word $\v$,  if $\delta(\widehat{\v^3})=0$ then $\bb(\wv)$
is empty.
\end{prop}

\begin{proof}  Consider a cyclically reduced, linear, nonpower word $\w$ such that  $\v=\w^k$ for some positive integer $k$.
By \lemc{linkedvp}, $\bb(\wv)$ is empty if an only if \(\bb(\widehat{\w})\)  is empty. Therefore, if  $k>1$, the conclusion follows from \propc{contar}. Hence, we can assume that $\v$ is nonpower.
Assume  $\bb(\wv)$ is not empty.
 Let $(i,j)$ be an element in $\bb(\wv)$. By Equation (\ref{eq delta}), one of the terms of $\delta(\widehat{\v^3})$ is  $3\cdot \sign(i,j)\widehat{\v_i^2\v_{i,j}}\otimes \wh{\v_{j,i} }$. Denote this term by $T$. We will show that $T$  does not cancel with other terms of $\delta(\widehat{\v^3})$.

Observe that in  Equation (\ref{eq delta}) there are two other terms corresponding to $(i,j)$, namely,
$3\cdot \sign(i,j)\widehat{\v_{i,j}}\otimes \wh{\v_j^{2}\v_{j,i} }$ and $3\cdot \sign(i,j)\widehat{\v_i\v_{i,j}}\otimes \wh{\v_j\v_{j,i} }$.
Clearly, none of these terms cancel with $T$.

By \propc{main}, if $(k,h) \in \bb(\wv)$ and $(k,h) \ne (i,j)$ then
$\widehat{\v
_k^{l}\v_{\hbox{\tiny $\!\!k\!,\!h$}}}  \ne \widehat{\v_i^2\v_{i,j}}$  for any positive integer $l$. Implies that the terms corresponding to $(k,h)$ do not cancel with $T$. Thus the proof is complete.\end{proof}

\section{Reformulation of Turaev's conjecture}\label{turaev}

Let  $\Sigma$ be an orientable surface with non-empty boundary.  One knows that there is a bijection between the following three sets:
\begin{numlist}
\item the set $\vvv$ of non-empty,  cyclic reduced words on a set of free generators of the fundamental group of $\Sigma$ and their inverses,
\item the set of non-trivial, free homotopy classes of closed, oriented curves on $\Sigma$, and
\item the set of conjugacy classes of $\pi_1(\Sigma)$.
\end{numlist}
We will make use of this correspondence by identifying these three sets. Thus an element in $\vvv$ will sometimes be considered as a conjugacy class of  $\pi_1(\Sigma)$ and other times, a free homotopy class  of closed, oriented curves.

Given a cyclic reduced word $\VV \in \vvv$, define the
\emph{self-intersection number  $\sss(\VV)$} to be minimum number of transversal self-intersection points of a representatives of  $\VV$ which self-intersect only in  transverse double points.

By \theoc{iso}, the combinatorially defined map $\delta$ (see \defc{cobr}) and the Turaev cobracket $\Delta$ (see \defc{Turaevcobracket}) coincide. By \theoc{iso}, \propc{simple} and \propc{contar} we have the following result.

\begin{proclama}{Main Theorem} Let $\Sigma$ be an oriented surface with boundary and let $\VV$ in $\vvv$. Then  $\VV$ contains a power of a simple curve if and only if $\Delta(\VV^3)=0$. Moreover, if $\VV$ is nonpower and $p$ is an integer larger than three then the Manhattan norm of $\Delta (\VV^p)$ equals $2p^2$ times the self-intersection number of $\VV$. In symbols,
$
\m(\Delta (\VV^p))=2p^2\sss(\VV).
$
\end{proclama}

\section{The Curve complex and the Mapping Class Group}\label{ccc}

Say that an element $\VV \in \vvv$ is \emph{simple} if contains an embedded curve. The following lemma will be used in the proof of \propc{permutation}.

 \start{lem}{eight} If $\UU$ and $\VV$  be two disjoint simple elements of  $\vvv$
 then there exists $P$ in  $\Sigma$ and   $\u$ and $\v$ in $\pi_1(\Sigma,P)$ such that $\u \in \UU$, $\v \in \VV$ and the following holds.
 \begin{numlist}
        \item Either $\wh{\u.\v}$ or $\wh{\u.\overline{\v}}$ has a representative with exactly one transverse double point.
        \item If $\UU\not=\VV$ and $\wh{\u.\v}$has a representative with exactly one transverse double point then $\wh{\u.\v}$ is nonpower and $\sss(\wh{\u.\v})=1$.
        \end{numlist}
\end{lem}

\begin{proof}  Let $\alpha$ and $\beta$ be disjoint closed simple curves that are representatives of $\UU$ and
$\VV$ respectively. Consider a connected component $C$ of $\Sigma\setminus(\alpha\cup\beta)$
including $\alpha\cup\beta$ in its boundary. Let $\gamma$ be a non self-intersecting  arc in $C$ from  a point $P$ in $\alpha$, to  a point $Q$ in $\beta$. Let $\u \in \pi_1(\Sigma,P)$ be the element  represented by $\alpha$ and let
$\v \in \pi_1(\Sigma,P)$ be the element represented by the curve that starts at $P$, runs along $\gamma$ to $Q$, then along $\beta$ and finally along $\gamma^{-1}$ from $Q$ to $P$.   It is not hard to see that $\u$ and $\v$ satisfy (1).

Now we prove (2). Since $\u.\v$ can be represented by a figure eight and $\UU$ and $\VV$ are non-trivial, $\sss(\wh{\u.\v})=1$. Consider an element $\WW \in \vvv$ such that $\wh{\u.\v}=\WW^p$,  for some positive integer $p$.  By definition of the cobracket,
$$
\Delta(\WW^p)=\Delta(\wh{\u.\v})=\pm(\UU \otimes \VV - \VV \otimes \UU) \ne 0.
$$
By \propc{power of simple}, $\WW$ is not simple, that is $\sss(\WW)>0$. By \cite[Theorem 3.9]{chas}, $ p^2\cdot s(\WW)\le
\sss(\WW^p)=\sss(\wh{\u.\v})=1$.
Then   $p=1$, and \(\wh{\u.\v}\) is
primitive.\end{proof}

 A permutation $\sigma$  on $\vvv$ \emph{ commutes with the power operation} if  for each $\VV$ in $\vvv$ and each $p$ in  $\Z$, $\sigma(\VV^p )=\sigma(\VV)^p$.
The linear map induced by $\sigma$ in the free $\Z$-module generated by $\vvv$ \emph{commutes with
Turaev's cobracket} if $\Delta(\sigma(\VV))=\sigma\otimes\sigma(\Delta(\VV))$.

\start{prop}{permutation} Let $\sigma$  be a  permutation on $\vvv$ which commutes
with the power operations  If the linear map induced by $\sigma$ on  $\mathbb Z[\vvv]$ commutes with
Turaev's cobracket then
  \begin{numlist}
           \item the map $\sigma$ preserves the subset of nonpower classes of a given self-intersection number. In particular, $\sigma$ maps  simple classes into simple classes, and
        \item  the restriction of $\sigma$ to the subset of simple classes maps pairs of classes with disjoint representatives into pairs of classes with disjoint representatives.
    \end{numlist}
 Therefore, $\sigma$ induces an automorphism on the curve complex $C(\Sigma)$ (see Appendix~\ref{cc}).
\end{prop}

\begin{proof}
Since $\sigma$ is bijective and commutes with the power operations,  the restriction of $\sigma$ to the subset of nonpower classes is a permutation.

 If  $\UU \in \vvv$ is  nonpower, $\sigma(\UU)$ is nonpower. Since $\sigma$ is a bijection,  by the Main Theorem,
$$
2\cdot4^2\!\cdot\!\sss\big(\sigma(\UU)\big)=\m(\Delta(\sigma(\UU)^4)=
\m(\sigma(\Delta(\UU^4)))=
\m(\Delta(\UU^4))=2\cdot4^2\!\cdot\!\sss(\UU).
$$
 Therefore, $\sss(\UU)=\sss(\sigma(\UU))$. This completes the proof of (1).

Finally, consider two simple classes $\UU$ and $\VV$ in $\vvv$ that can be represented by disjoint closed curves.  If $\UU=\VV$, then the conclusion is immediate. Thus we can assume $\UU \ne \VV$.
By \lemc{eight}, there exist  $\u$ and $\v$ in $\pi_1(\Sigma)$ satisfying \lemc{eight}(1) and (2).
 Suppose first that $\wh{\u\cdot\v}$ can be represented by a figure eight. By (1), $\sss(\sigma(\wh{\u\cdot\v}))=1$. Then there exists two elements $\x$ and $\y$ in $\pi_1(\Sigma)$ such that $\sigma(\wh{\u\cdot\v})=\wh{\x\cdot\y}$. Moreover, $\wh{\x}$ and $\wh{\y}$ have simple disjoint representatives. Hence,
$$\pm(\wh{\x}\otimes \wh{\y}-\wh{\y}\otimes \wh{\x})=\Delta(\sigma(\wh{\u\cdot \v}))= \sigma\otimes\sigma(\Delta(\wh{\u\cdot \v})) =\pm(\sigma( \wh{\u}) \otimes \sigma(\wh{\v})-\sigma(\wh{\v}) \otimes \sigma(\wh{\u}) )$$
Then  the set  $\{\wh{\x},\wh{\y}\}$ is equal to the set $\{\sigma(\wh{\u}),\sigma(\wh{\v})\}$. This implies that $\sigma(\wh{\u})$ and $\sigma(\wh{\v})$ are disjoint.

If $\wh{\u \cdot \overline{\v}}$ is represented by a figure eight, then by the arguments above we can show that $\sigma(\u)$ and $\sigma(\overline{\v})$ are disjoint. Since  $\sigma(\overline{\v})=\overline{\sigma(\v)}$, $\sigma(\u)$ and $\sigma(\v)$ are disjoint.
\end{proof}

We  denote by $\Sigma_{g,b}$ an oriented surface with genus $g$ and $b$ boundary components.

\start{theo}{comples} Let $\sigma$ be a permutation  on the set $\pi^\ast$ of free homotopy  classes of closed
curves on an oriented surface with boundary which is different from $\Sigma_{1,2}$, $\Sigma_{0,4}$ and $\Sigma_{1,1}$.  Suppose that $\sigma$ commutes with the power operation and that the linear function induced by $\sigma$ in $\Z[\pi^*]$ commutes with Turaev's cobracket.
Then $\sigma$ induces a permutation $\widetilde{\sigma}$ on the set of unoriented simple classes.
The permutation  $\widetilde{\sigma}$  is induced by a unique element of the mapping
class group. \end{theo}
\begin{proof}
By \propc{permutation}(1), since $\sigma(\overline{\UU})=\overline{\sigma(\UU)}$, $\sigma$ induces a permutation $\widetilde{\sigma}$ on the set of unoriented simple closed curves.

By \propc{permutation}(2), the restriction of $\widetilde{\sigma}$ to the set of unoriented simple closed curves preserves disjointness.

Thus, by \theoc{ikl}, $\widetilde{\sigma}$ is induced by an element of the Mapping Class group.
\end{proof}

This result ``supports'' Ivanov's statement in \cite{iv2}:

\begin{proclama}{Metaconjecture} ``Every object naturally associated with a surface $S$ and having a sufficiently rich
structure has $\mathrm{Mod}(S)$ as its group of automorphisms. Moreover, this can be proved by a reduction
theorem about the automorphisms of $\mathrm{C}(S)$.''
\end{proclama}

In this sense, the Turaev Lie cobracket combined with the power maps, have a ``sufficiently rich" structure.

\appendix

\section{Turaev's Lie coalgebra of curves on a surface}\label{curves}

Consider a vector space  $\W$. Define two linear maps $\Smap{\om}{\W\otimes \W\otimes \W}$ and $\Smap{\sss}{\W\otimes \W}$
by the formulae $\om(u \otimes v\otimes w) =w \otimes u \otimes v$ and $\sss(v\otimes
w)=w\otimes v$ for each triple of elements $u$, $v$ and $w$ in $\W$.
A {\em Lie cobracket} on $\W$ is a linear map $\map{\cib}{\W}{\W\otimes
\W}$ such that $\sss \circ\cib=-\cib$ ({\em co-skew symmetry}) and $(\id+\om+\om^2)(\id
\otimes \cib)\cib=0$ ({\em co-Jacobi identity}). If $\W$  is a vector space and $\cib$ is a cobracket on $\W$, $(\W,\cib)$ is a {\em Lie coalgebra}.

We recall the definition of Turaev's cobracket \cite{T}.

\start{defi}{Turaevcobracket} Let $\Sigma$ be an oriented surface. Choose a free set of generators for the fundamental group of $\Sigma$. The set of non-trivial free homotopy classes of curves on $\Sigma$ is in bijective correspondence with the set of cyclic  words (as defined in  Section \ref{tur}) in the generators and their inverses. We will identify cyclic words with free homotopy classes.
Thus $\V$  is the vector space generated by the set of free homotopy classes of
non-trivial, oriented, closed curves on $\Sigma$.

Let $\gamma$ be a oriented, closed curve on $\Sigma$. Let $[\gamma]$ denote the free homotopy class of $\gamma$ if $\gamma$ is a non-contractible loop and $[\gamma]=0$ otherwise.
Consider a free homotopy class $[\beta]$ such that all self-intersection points of $\beta$ are  transversal double points.  Denote by $\mathcal{I}$ the set of self-intersection points of $\beta$. Each $P \in \mathcal{I}$ determines two loops based at $P$ so that $\beta$ can be obtained as the loop product of these two loops (forgetting the basepoint). Order these two loops so that the orientation given by the branch of the first, followed by the branch of the second equals the orientation of the surface. Denote the ordered pair of loops by   $(\beta_P^1,\beta_P^2)$.
The Turaev cobracket of $[\beta]$, $\Delta([\beta])$ is defined by the following formula:
\[
\Delta([\beta])=\sum_{P\in \mathcal{I}}
            \Big([\beta_P^1]\otimes [\beta_P^2]
            -[\beta_P^2]\Big)\otimes[\beta_P^1]
\]

\end{defi}
The next result follows straightforwardly from the definition of the cobracket.
\start{prop}{power of simple} If a free homotopy class of curves $x$ can be represented by a power of a simple curve then $\Delta(x)=0$.
\end{prop}

\section{Examples }\label{review c}

\begin{figure}[htbp]
\begin{center}
	\includegraphics[scale=0.27]{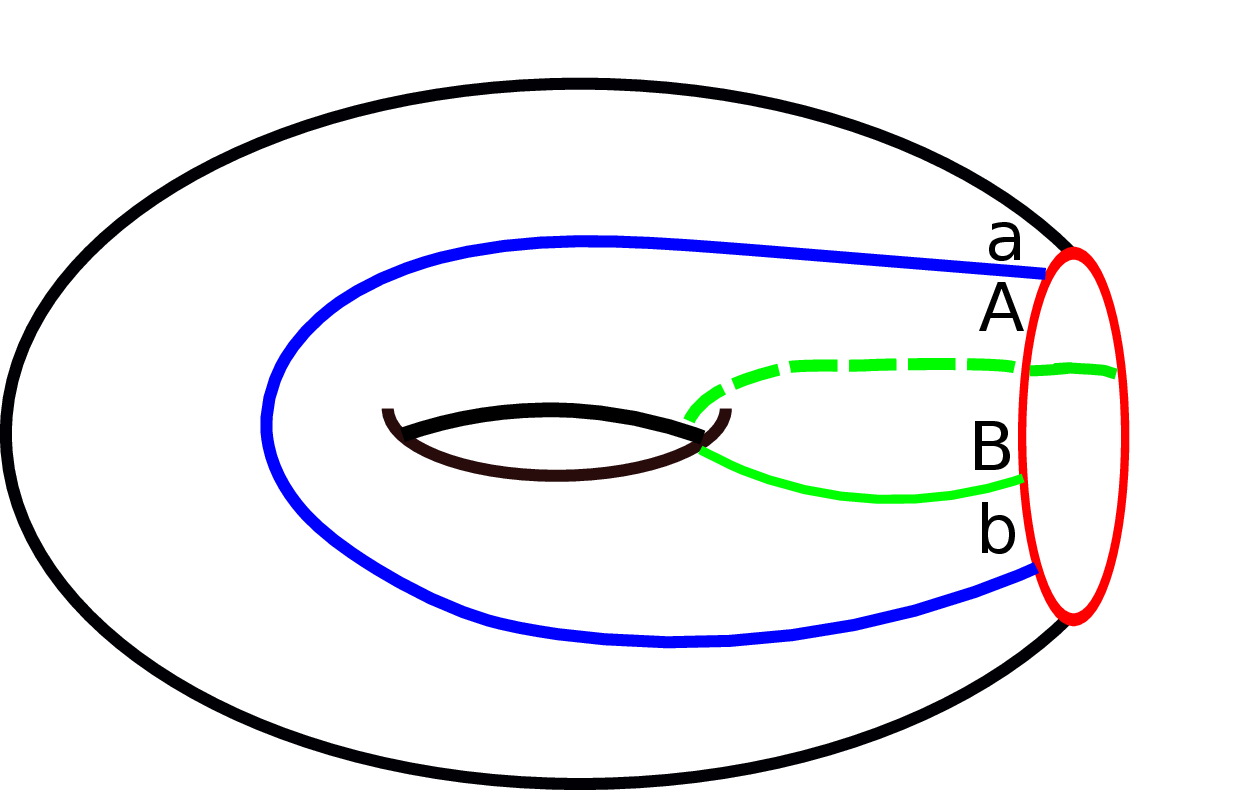}
\includegraphics[scale=0.27]{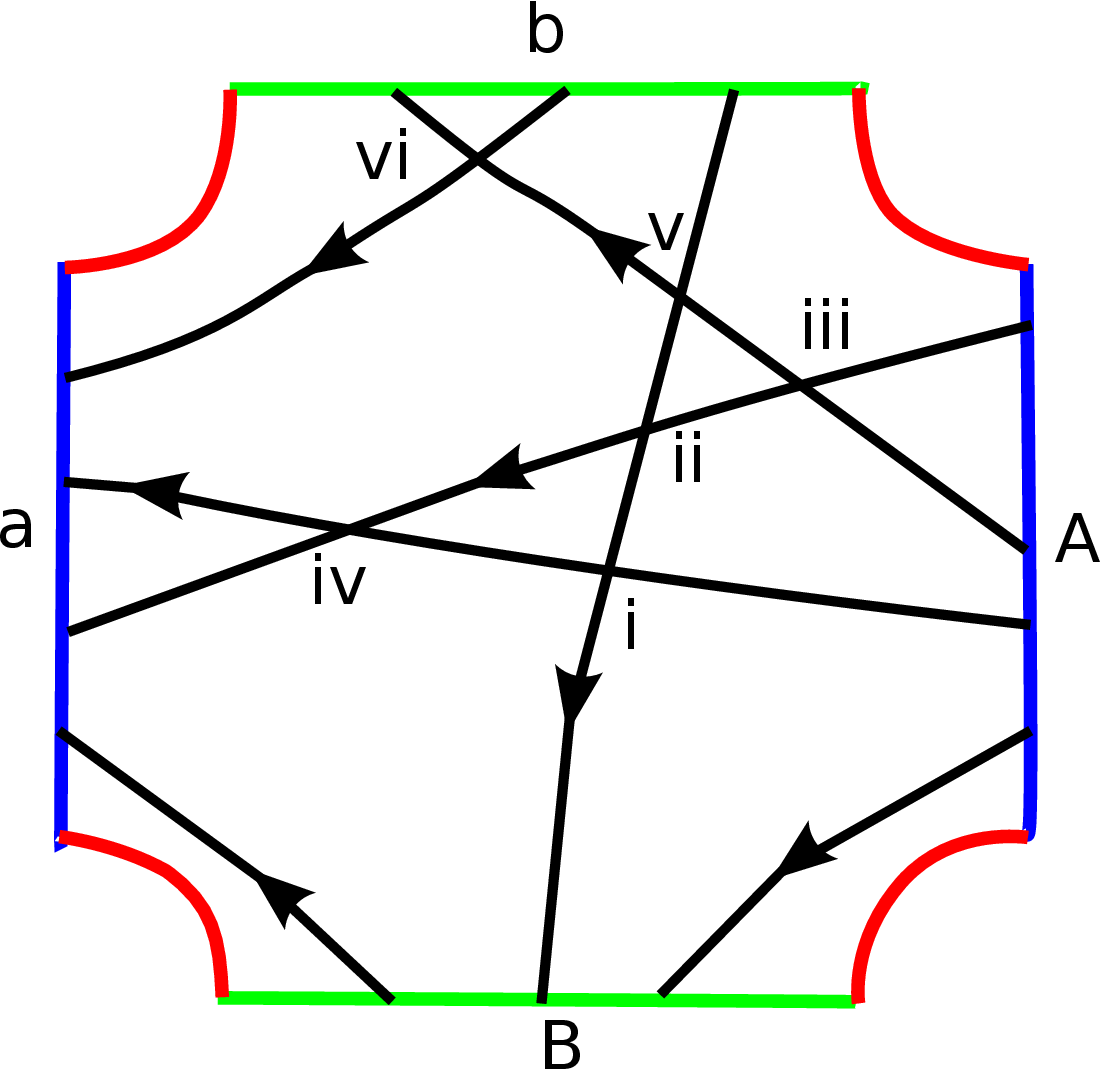}
\caption{The arcs in the punctured torus (left) and the  curve $aaabaBB$ (right) }
\label{torusCurve}
   \end{center}
 \end{figure}

Consider an orientable surface  with non-empty boundary and negative Euler characteristic.  Choose a maximal set of disjoint arcs each starting and ending in the boundary, such that the surface minus the union of the arcs is connected (see Figure~\ref{torusCurve}, left). Note that this connectivity plus maximality implies that no two of the arcs are homotopic keeping endpoints in the boundary (i.e., rel. boundary).

Label one side of each arc with a letter $x$, and the other side with the letter $\bar{x}$.  The choice of arcs determines a minimal set of generators of the fundamental group of the surface: Choose a basepoint $P$ in the surface in the complement of the chosen arcs. A representative of the generator labeled by $x$ is a curve that starts at $P$ crosses the arc labeled  $x$ from the side labeled by $x$ to the side labeled by $\bar{x}$ and goes back to $P$ without crossing any other arc. By cutting the surface along these arcs, a polygon is obtained with four times the number of sides as the number of generators (see Figure~\ref{torusCurve}, right). Alternating edges are labeled.    By choosing an edge and reading the labels of the edges in cyclic order a linear order of the edge labels is obtained. In the example of Figure~\ref{torusCurve} the order  is $a,b,A,B$. (To ease the notation, we write $A$ instead of $\overline{a}$ and $B$ instead $\overline{b}$

This linear order of the alphabet (called the surface word) determines the structure of  self-intersection of the free homotopy classes given in terms of the above generators of the fundamental group. For instance, in the torus with one boundary component of our example,  $\widehat{aaabaBB}$ since the  curve word contains the pair subwords $aa$ and $BB$, any representative of that class will have an arc from the edge labeled by $a$ to the edge labeled by $A$  and another arc from the edge labeled by $b$ to the edge labeled by $B$. These two arc must intersect. Thus the subwords $aa$ and $BB$ "imply" a self-intersection point. This is a key idea.

The linear order also determines an orientation of the surface, and so, a sign at each intersection point of a minimal representative.

\start{ex}{exa} Half of the linked pairs of the word $\v=\widehat{aaabaBB}$, types and cobracket terms  are listed in Table~\ref{table}. (The rest of the linked pairs are obtained by swapping $\p$ and $\q$.  A representative of this class with labeled linked pairs is displayed in Figure~\ref{torusCurve}.
\end{ex}

\begin{table}
\centering
\begin{tabular}{|c|c|c|c|c|c|c|c|}
\hline
\rule[-1ex]{0pt}{2.5ex}  & $\p$&$\q$ & (i,j) & $Y$ &Type & Cobracket Term &  $\sign(\p,\q)$\\
\hline
\rule[-1ex]{0pt}{2.5ex} i &  {\it aa}& \textit{BB} & (1,6) & -&(1) & $\widehat{aabaB} \otimes \widehat{Ba}$ &-1 \\
%\hline
\rule[-1ex]{0pt}{2.5ex} ii &\textit{aa} &\textit{BB} & (2,6) & -&(1) &$\widehat{abaB} \otimes \widehat{Baa}$ & -1\\
%\hline
\rule[-1ex]{0pt}{2.5ex} iii& \textit{aab} &\textit{Baa} & (3,1) & $a$&(2)  & $\widehat{baBBa} \otimes \widehat{aa}$& -1\\
%\hline
\rule[-1ex]{0pt}{2.5ex} iv & \textit{Baaa}& \textit{aaab} & (2,3) & $aa$ &(2) &$\widehat{a} \otimes \widehat{baBBaa}$ & 1\\
%\hline
\rule[-1ex]{0pt}{2.5ex} v & \textit{aba} & \textit{aBB} & (4,5) & $b$&(3) & $\widehat{a} \otimes \widehat{Baaa}$& -1 \\
%\hline
\rule[-1ex]{0pt}{2.5ex} vi & \textit{aba} & \textit{BBa} & (4,6) & $b$&(3) &$ \widehat{aB} \otimes \widehat{aaa}$ & -1 \\
\hline
\end{tabular}
\caption{
Six of the twelve linked pairs of $\v=aaabaBB$ and the corresponding cobracket terms. The labels $i, ii, \dots, vi$ refer to Figure~\ref{torusCurve}. }\label{table}
\end{table}

\begin{figure}[hp]
\begin{center}
	\includegraphics[scale=0.7]{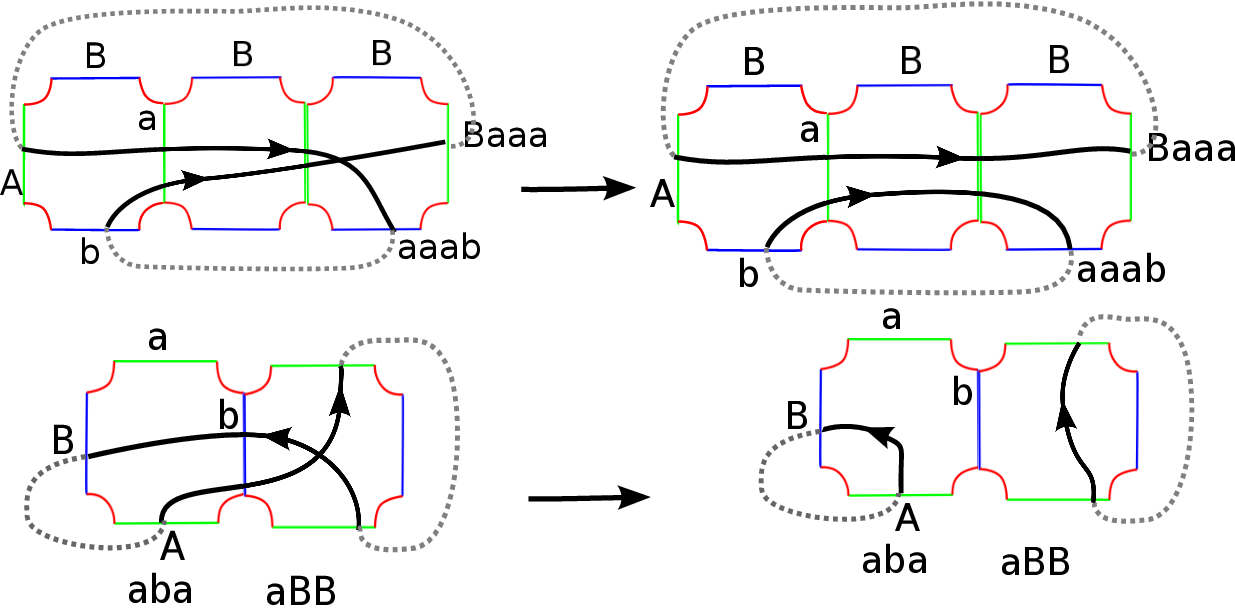}
\caption{The horizontal arrows perform the cobracket of the linked pairs $(aaab,Baaa)$ and $(aba,aBB)$}
\label{lpexample}
   \end{center}
 \end{figure}

\start{ex}{exa1} In Figure~\ref{lpexample}, a representation of the cobracket terms corresponding to $iv$ and $v$ in Figure~\ref{torusCurve} are exhibited. The pair $iv$, determines the ``cuts'' $aa|a|baBB$. To obtain the cobracket term, one considers the words $a$ and $aabaBB$ and makes them cyclic. The pair $v$ determines the ``cuts'' $aaab|a|BB$. Since this is a pair of type $(3)$. after ``cutting'' one should reduce the superfluous part of one of the words. The obtained words are $aaaB$ and $a$.
\end{ex}

\begin{figure}[hpb]
\begin{center}
\includegraphics[scale=0.5]{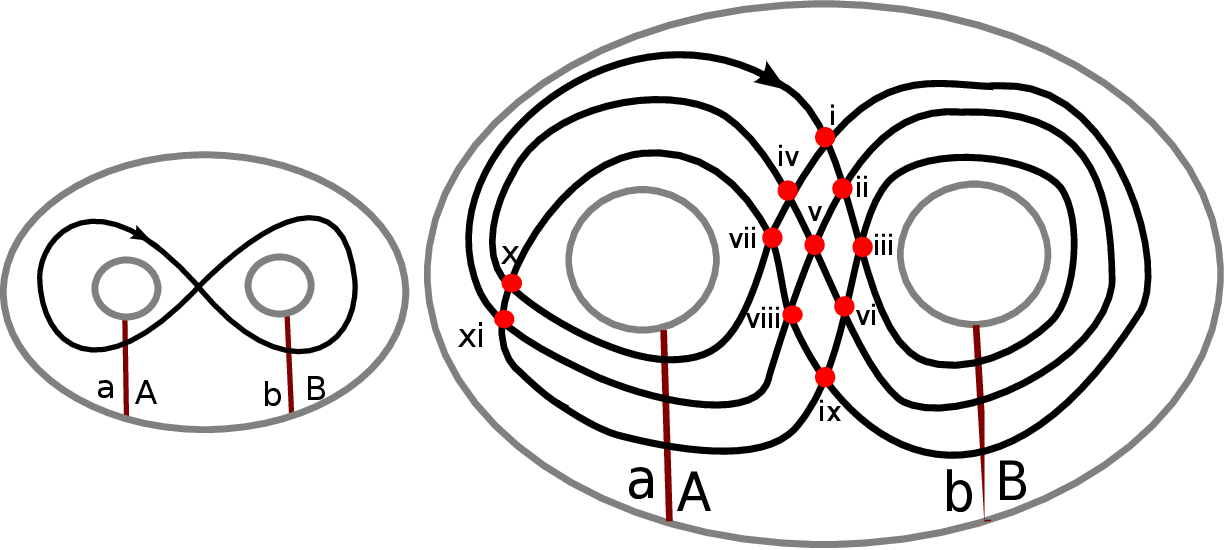}
\caption{Linked pairs of the third power word $bAbAbA$ (right) and the root word $bA$ (left)}
\label{lpexample1}
   \end{center}
 \end{figure}

\start{ex}{exa2} Six of the eighteen  linked pairs of the power word  $\v=bAbAbA$   are listed in Table~\ref{table1}.  A representative of this class is displayed in Figure~\ref{lpexample1}. Note that in the geometric version of the cobracket the terms corresponding to the points labeled  $x$ and $xi$   in Figure~\ref{lpexample1} cancel. The linked pairs do not ``see'' these two points.
\end{ex}

\begin{table}
\centering
\begin{tabular}{|c|c|c|c|c|c|c|c|}
\hline
\rule[-1ex]{0pt}{2.5ex}  & $\p$&$\q$ & (i,j) & $Y$ &Type & Cobracket Term&  $\sign(\p,\q)$\\
\hline
\rule[-1ex]{0pt}{2.5ex} i &  \textit{Ab}& {\it bA}  & (0,3) & -&(1)&
$\widehat{bab}\otimes \widehat{Aba}$ &-1\\
\hline
\rule[-1ex]{0pt}{2.5ex} ii &\textit{Ab}& {\it bA} & (0,5) & -&(1)&
$\widehat{bAbAb}\otimes \widehat{A}$ &-1\\
\hline
\rule[-1ex]{0pt}{2.5ex} iii& \textit{Ab}& {\it bA} & (0,1) & -&(1) &
$\widehat{b}\otimes \widehat{AbAbA}$&-1\\
\hline
\rule[-1ex]{0pt}{2.5ex} iv & \textit{Ab}& {\it bA} & (4,3) & -&(1)&
$\widehat{bAbAb}\otimes \widehat{A}$ &-1\\
\hline
\rule[-1ex]{0pt}{2.5ex} v & \textit{Ab}& {\it bA} & (4,5) & -&(1)&
$\widehat{b}\otimes \widehat{AbAbA}$ &-1\\
\hline
\rule[-1ex]{0pt}{2.5ex} vi & \textit{Ab}& {\it bA} & (4,1) & -&(1)&
$\widehat{bAb}\otimes \widehat{AbA}$ &-1\\
\hline
\end{tabular}
\caption{Six of the eighteen linked pairs of $\v=bAbAbA$ and the correspoding cobracket terms. The labels $i, ii, \dots vi$ refer to Figure~\ref{lpexample1}.}\label{table1}
\end{table}

\section{The curve complex}\label{cc}

Let $\Sigma$ be a compact oriented surface.  By $\Sigma_{g,b}$ we denote an oriented surface with genus $g$
and $b$ boundary components. If $\Sigma$ is a surface, we denote by $\mathcal{MCG}(\Sigma)$ the \emph
{mapping class group of $\Sigma$}, that is, the set of homotopy classes of orientation preserving homeomorphisms of $
\Sigma$.  We study automorphisms of the Turaev Lie coalgebra that are related to the mapping class group.

Now we recall the curve complex,  defined by Harvey in \cite{wh}. The \emph{curve complex   $\mathrm{C}(\Sigma)$ of $\Sigma$}
 is the simplicial complex  whose vertices are  isotopy classes of unoriented  simple closed
curves on $\Sigma$ which are neither null-homotopic nor homotopic to a boundary component. If $\Sigma \ne
\Sigma_{0,4} $ and $\Sigma \ne \Sigma_{1,1}$ then a set of $k+1$ vertices of the curve complex is the $0$-
skeleton of a $k$-simplex if the corresponding  minimal intersection number of all pairs of vertices is zero, that is, if
every pair of vertices have disjoint representatives.

For $\Sigma_{0,4}$ and $\Sigma_{1,1}$ two vertices are connected by an edge when the curves they represent
have minimal intersection (2 in
the case of $\Sigma_{0,4}$ , and 1 in the case of $\Sigma_{1,1}$).  If $b \le 3$ the complex associated with $\Sigma_{0,b}$ is empty.

The following isomorphism is a theorem of Ivanov
\cite{iv} for the case of genus at least two.  Korkmaz \cite{ko} proved the result for genus at most one and Luo \cite{fl} gave
another proof  that covers all possible genera. The mapping class group of a surface $\Sigma$ is denoted by $\mathcal{MCG}(\Sigma)$. Our statement below is based on the formulation of Minsky \cite{minsky}.

\start{theo}{ikl}(Ivanov-Korkmaz-Luo) The natural map $\map{h}{\mathrm{\mathcal{MCG}}(\Sigma)}{\mathrm{Aut C}(\Sigma)}$ is an isomorphism in all cases except for $\Sigma_{1,2}$ where it is injective with index $2$ image.
 \end{theo}


\begin{thebibliography}{99}

\bibitem{Ar} C. Arettines, A combinatorial algorithm for visualizing representatives with minimal self-intersection
arXiv:1101.5658, to appear in Journal of Knot Theory and Ramifications.

\bibitem{BS} J.\ Birman and C.\ Series, An algorithm for simple curves on surfaces, {\em J. London Math. Soc.} (2), {\bf 29} (1984), 331-342.


\bibitem{Cahn} Cahn, P.,
\textit{ A Generalization of the Turaev Cobracket and the Minimal Self-Intersection Number of a Curve on a Surface.}
\href{http://arxiv.org/abs/math/1004.0532}{{\tt arXiv: 1004.0532}}




\bibitem{CK} Chas, M., Krongold, F.,
\textit{ An algebraic  characterization of simple closed curves on
surfaces with boundary},  Journal of Topology and Analysis (JTA), Volume: 2, Issue: 3(2010) pp. 395-417
World Scientic Publishing Company.
%\href{http://front.math.ucdavis.edu/0801.3944}{{\tt arXiv:0801.3944
%[math.GT]}}.



\bibitem{chas} Chas, M., \textit{Combinatorial Lie bialgebras of curves on surfaces},
Topology \textbf{43},   (2004), 543-568.
\href{http://arxiv.org/abs/math/0105178v2}{{\tt arXiv: 0105178v2
[math.GT]}}

\bibitem{cl} M.\ Cohen and M.\ Lustig, Paths of geodesics and geometric intersection numbers I,  {\em Combinatorial Group Theory and Topology}, Alta, Utah, 1984,  Ann. of Math. Studies  {\bf 111}, Princeton Univ. Press, Princeton,  (1987), 479-500.

\bibitem{chi1} D. Chillingworth,  Winding numbers on surfaces. I, Math. Ann. 196 (1972), 218-249

\bibitem{chi2} D. Chillingworth,  Winding numbers on surfaces. II, Math. Ann. 199 (1972), 131-153.




\bibitem{wh} W.~J. Harvey, {Geometric structure of surface mapping class groups,} {Homological group theory (Proc. Sympos., Durham, 1977)}, 255--269. Cambridge Univ. Press, Cambridge, 1979.

\bibitem{hs} J. Hass, P. Scott, \textit{Intersections of curves on surfaces, } Israel Journal of mathematics \textbf{51}, No.1-2, (1985).

\bibitem{hs1}	J. Hass and P. Scott, Shortening curves on surfaces, Topology, Vol 33, \textbf{1} 25-43, (1994).

\bibitem{iv}    N. Ivanov, {Automorphisms of complexes of curves and of {T}eichm\"uller spaces,} {Progress in knot theory and related topics}, 113--120. Hermann, Paris, 1997.
\bibitem{iv2} N. Ivanov, {Fifteen problems about the Mapping Class group,} Problems on mapping class groups and related topics, (editor Benson Farb), 71-80,  Providence, R.I., American Mathematical Society, (2006).


\bibitem{jaco} W. Jaco,  Heegaard splittings and splitting homomorphisms. Trans. Amer. Math. Soc. 144 (1969) 365-379.


\bibitem{ko} M. Korkmaz, {Automorphisms of complexes of curves on
punctured spheres and on punctured tori.} Topology Appl., {\bf  95}(2):85--111, (1999).
\bibitem{ald} A. Le Donne,  On Lie bialgebras of loops on orientable surfaces., J. Knot Theory Ramifications, 17, 351 (2008).


\bibitem{ledonne} A. Le Donne, On Lie bialgebras of loops on orientable surfaces.  J. Knot Theory Ramifications {\bf 17} (2008), no. 3, 351--359.

\bibitem{l} M. Lustig, Paths of geodesics and geometric intersection numbers II, Combinatorial
group theory and topology (Alta, Utah, 1984), Ann. of Math. Stud. 111, Princeton Univ.
Press, Princeton, NJ, 1987, pp. 501-543.

%\bibitem{luo1} F. Luo, {On non-separating simple closed curves in a compact surface,} Topology {\bf 36}, No 2,pp 381-410, (1997).
\bibitem{fl}       F. Luo, {Automorphisms of the complex of curves.} Topology {\bf 39}(2), 283--298, (2000).  \href{http://arxiv.org/abs/math/9904020}{{\tt arXiv: 9904020 [math.GT]}}
\bibitem{minsky} Y. Minksy, {Curve complexes, surfaces and 3-manifolds}, International Congress of
Mathematicians, Vol.{\bf II}, 1001-1033, Eur. Math. Soc., Zurich, (2006).

\bibitem{poincare} H. Poincar\`e , Analysis situs, J. Ecole Polytechn. (2), 1 (1895), 1-121.
%\bibitem{tu}  V. Turaev,


\bibitem{stallings} J. Stallings, How not to prove the Poincar\`{e} conjecture, Topology Seminar, Wisconsin, 1965. Edited by R. H. Bing and R. J. Bean. Annals of Mathematics Studies, No. 60. Princeton University Press, Princeton, N.J. 1966

\bibitem{T} V,  Turaev,
\textit{Skein quantization of Poisson algebras of loops on
surfaces}, Ann. Sci. Ecole Norm. Sup. (4) \textbf{ 24}, No. 6,
 (1991), 635-704 .
\bibitem{tv}  V. Turaev, O. Viro, Intersection of loops in two-dimensional manifolds. II. Free loops. Mat. Sbornik 121:3 (1983) 359-369 (Russian); English translation in Soviet Math. Sbornik.
\bibitem{z1} H. Zieschang,  Algorithmen fur einfache Kurven auf Fl\"{a}chen, Math. Scand. 17 (1965),
17-40. 21.
\bibitem{z2} H. Zieschang,  Algorithmen fur einfache Kurven auf Fl\"{a}chen. II, Math. Scand. 25 (1969),
49-58.
\end{thebibliography}
\end{document}